# Test for tail index change in stationary time series with Pareto-type marginal distribution

MOOSUP KIM[*] and SANGYEOL LEE[**]

*Department of Statistics, Seoul National University, Seoul, 151-742, Korea.*
*E-mail: [*]mooseob@hanmail.net; [**]sylee@stats.snu.ac.kr*

The tail index, indicating the degree of fatness of the tail distribution, is an important component of extreme value theory since it dominates the asymptotic distribution of extreme values such as the sample maximum. In this paper, we consider the problem of testing for a change in the tail index of time series data. As a test, we employ the cusum test and investigate its null limiting distribution. Further, we derive the null limiting distribution of the cusum test based on the residuals from autoregressive models. Simulation results are provided for illustration.

## 1. Introduction

The parameter change problem in time series models has attracted much attention from researchers since time series often experience changes in the underlying models due to the changes of monetary policies and critical social events. Since Page (1955), a large number of studies have been devoted to the theory and its applications on the change point analysis in various fields. For a general review of the change point analysis, we refer to Csörgő and Horváth (1997). Among the existing methods, the cusum test has long been popular for its ease of usage in actual practice. Compared to the likelihood method, the cusum method has an advantage that the null limiting distribution is free from the underlying distribution. For relevant references, we employ the following articles and the literatures cited therein: Brown, Durbin and Evans (1975); Tang and MacNeil (1993); Inclán and Tiao (1994); Lee and Park (2001); Lee, Ha, Na and Na (2003); Berkes, Horváth and Kokoszka (2004); Lee and Na (2005); and Lee, Nishiyama and Yoshida (2006).

Statistical modeling and analysis for extremal phenomena is very crucial in that the potential risk of disasters and panic events such as floods, large earthquakes and stock market crashes can be determined a priori, thereby allowing them to be adequately managed or prevented. On the other hand, it is well known that structural changes in the underlying models can lead to false conclusions as frequently observed in the financial time series analysis. Motivated by this, we are led to study the change point test in the extreme value theory. In particular, we focus on the change point test for the tail index of distributions since the tail index represents the degree of fatness of distributions and determines the shape of the asymptotic distribution of extreme values such as the sample maximum.

The estimation problem for the tail index has been a core issue for several decades in statistics, finance, reliability and teletraffic engineering. Among others, Hill's estimator (cf. Hill (1975)) has been playing an important role in this context and its properties are well developed in various respects. In fact, one of the crucial issues with regard to Hill's estimator is how to select the tail sample fraction since a large sample fraction yields more bias in Hill's estimator. For a general review of the extreme value theory and its statistical applications, we refer to Leadbetter, Lindgren and Rootzén (1983); Embrechts, Küppelberg and Mikosch (1999); and Reiss and Thomas (2001). Further, for the references relevant to Hill's estimator, we refer to Hall (1982); Mason (1982); Csörgő, Deheuvels and Mason (1985); Smith (1987); Hsing (1991); and Resnick and Stărică (1995, 1997a, 1997b).

For performing the change point test, here we employ the cusum test based on Hill's estimator. In particular, we concentrate on the cusum test for the $\beta$-mixing process that includes a broad class of stationary processes such as the autoregressive (AR), generalized autoregressive conditional heteroscedastic (GARCH) and threshold AR (TAR) processes (cf. Doukhan (1994)). Further, for the AR process, we consider the residual-based cusum test since the tail index of the error distribution coincides with that of the AR process itself (cf. Datta and McCormick (1998)). In general, the residual-based test is more stable compared to those based on the observations themselves due to the elimination of the correlation effects (see Resnick and Stărică (1997a) for Hill's estimator case). In order to construct the cusum test, we consider the tail sequential process in light of the testing procedure proposed by Lee, Ha, Na and Na (2003), and verify that the tail sequential process converges weakly to a Brownian motion under the null hypothesis under which the tail index is assumed to remain as a constant. In this case, the asymptotic null distribution appears to be the sup of a Brownian bridge.

This paper is organized as follows: In Section 2, we introduce the cusum test and derive its asymptotic null distribution. Further, we discuss the estimation procedure for the location of change points. In Section 3, we perform a simulation study to evaluate our tests. In Section 4, we provide the proofs for the results presented in Section 2.

## 2. Main results

### 2.1. Cusum test in $\beta$-mixing processes

In this subsection, we consider the problem of testing for the change of the tail index in a class of $\beta$-mixing processes. In what follows, we assume that all r.v.'s are defined on the probability space $(\Omega, \mathcal{F}, P)$. Let $\{X_i\}$ be a sequence of non-negative r.v.'s. Suppose that one wishes to test the following hypotheses:

$\mathcal{H}_0 : \{X_i\}$ is stationary and the tail index does not vary over $X_1, \ldots, X_n$ vs.

$\mathcal{H}_1 :$ not $\mathcal{H}_0$.

Under $\mathcal{H}_0$, we assume that $\{X_i\}$ satisfies the $\beta$-mixing condition:

$$\beta(l) = \sup_m \mathrm{E}\left\{\sup_{A \in \mathcal{F}_{m+l+1}^\infty} |P(A|\mathcal{F}_1^m) - P(A)|\right\} \to 0 \quad \text{as } l \to \infty, \tag{2.1}$$

where $\mathcal{F}_l^m = \sigma\{X_i : i = l, \ldots, m\}$ and $\mathcal{F}_m^\infty = \sigma\{X_i : i = m, m+1, \ldots\}$. Further, we assume that the common marginal distribution $F$ of $\{X_i\}$ has the tail index $\alpha^{-1} > 0$, namely,

$$\lim_{x \to \infty} \frac{\bar{F}(\lambda x)}{\bar{F}(x)} = \lambda^{-\alpha} \quad \text{for each } \lambda > 0, \tag{2.2}$$

where $\bar{F} = 1 - F$. In this case, $\bar{F}$ is said to be regularly varying at $\infty$ with the exponent $-\alpha$ (abbreviated as $\bar{F} \in RV_{-\alpha}$). According to Theorem 1.6.2 of Leadbetter *et al.* (1983), $\bar{F}$ is regularly varying at $\infty$ if and only if $F$ lies in the domain of attraction of the Fréchet distribution. Owing to (2.2), we can express $\bar{F}(x) = x^{-\alpha} l(x)$, where $l$ is slowly varying at $\infty$, namely,

$$\lim_{x \to \infty} \frac{l(\lambda x)}{l(x)} = 1 \quad \text{for every } \lambda > 0. \tag{2.3}$$

For performing a test for $\mathcal{H}_0$ and $\mathcal{H}_1$, we employ the cusum test based on Hill's estimator:

$$\frac{1}{k} \sum_{i=1}^n (\log X_i - \log X_{(k+1)})_+, \tag{2.4}$$

where $x_+ = \max\{x, 0\}$, $X_{(j)}, j = 1, \ldots, n$, denotes the $j$th largest r.v. in $X_1, \ldots, X_n$, and $k$ is a positive integer much less than $n$. According to Lee, Ha, Na and Na (2003), the cusum test is constructed based on the following tail sequential process:

$$M_n(t) := \frac{1}{\sqrt{k}} \sum_{i=1}^{[nt]} \{\varphi(\log X_i - \log b(n/k)) - \mathrm{E}\varphi(\log X_i - \log b(n/k))\}, \quad 0 \leq t \leq 1, \tag{2.5}$$

where $\varphi$ is a real-valued function with $\varphi(x) = 0$ for every $x < 0$, $b(x) = \inf\{y : F(y) \geq 1 - x^{-1}\}$, and $k = k_n$ is a sequence of positive integers satisfying

$$k \to \infty \quad \text{and} \quad k = \mathrm{o}(n) \tag{2.6}$$

as $n \to \infty$. In this study, we particularly concentrate on the two cases:

$$\varphi_1(x) := I(x > 0) \quad \text{and} \quad \varphi_2(x) := x_+.$$

In order to investigate the limiting behavior of $M_n$, we assume that the following regularity conditions hold:

**(A1)** There exists a sequence $\{r_n\}$ of positive integers such that

$$\lim_{n \to \infty} \frac{n}{r_n} \beta([\varepsilon r_n]) = 0 \qquad \text{for every } \varepsilon > 0, \tag{2.7}$$

and

$$r_n^2 = \mathrm{o}(k). \tag{2.8}$$

**(A2)** There exist $\kappa(x) = K \int_1^x t^{\gamma-1} \, \mathrm{d}t$, where $K \in \mathbb{R}$, and a positive measurable function $g \in RV_\gamma$ with $\gamma \leq 0$, such that for all $\lambda > 0$,

$$\lim_{x \to \infty} \frac{l(\lambda x)/l(x) - 1}{g(x)} = \kappa(\lambda).$$

Further, $\sqrt{k} g(b(n/k))$ converges to a real number $A$ as $n \to \infty$.

**(A3)** There exist non-negative numbers $\chi$ and $\omega$ such that for every $0 < \varepsilon < 1$,

$$\chi = \lim_{n \to \infty} \frac{2\alpha^2 n}{[\varepsilon r_n] k} \times \sum_{1 \leq i < j \leq [\varepsilon r_n]} \mathrm{Cov}\{(\log X_i - \log b(n/k))_+, (\log X_j - \log b(n/k))_+\} \tag{2.9}$$

and

$$\omega = \lim_{n \to \infty} \frac{2n}{[\varepsilon r_n] k} \sum_{1 \leq i < j \leq [\varepsilon r_n]} \mathrm{Cov}\{I(X_i > b(n/k)), I(X_j > b(n/k))\}. \tag{2.10}$$

Condition **(A2)** is referred to as the second-order regularly varying condition and $\gamma$ is called the second-order regularly varying parameter. This condition plays a crucial role in the derivation of the asymptotic properties for tail index estimators. For the details concerning **(A2)**, readers are referred to Bingham *et al.* (1987) and Goldie and Smith (1987). In fact, it can be easily seen that Condition **(A3)** is satisfied for a large class of short memory processes.

According to our analysis, it is revealed that under the regularity conditions $M_n(t)$ divided by a constant converges weakly to a standard Brownian motion and, subsequently, the following random sequence:

$$T_n^\circ(\varphi) := \frac{1}{\sqrt{k}} \max_{1 \leq l \leq n} \left| \sum_{i=1}^l \varphi(\log X_i - \log b(n/k)) - \frac{l}{n} \sum_{i=1}^n \varphi(\log X_i - \log b(n/k)) \right| \tag{2.11}$$

converges weakly to the sup of a Brownian bridge multiplied by a constant. Since $b(n/k)$ is unknown, we replace it with $X_{(k)}$ and finally employ the cusum test statistic:

$$T_n(\varphi) := \frac{1}{\sqrt{k}} \max_{1 \leq l \leq n} \left| \sum_{i=1}^{l} \varphi(\log X_i - \log X_{(k)}) - \frac{l}{n} \sum_{i=1}^{n} \varphi(\log X_i - \log X_{(k)}) \right|. \tag{2.12}$$

In particular, we can express

$$T_n(\varphi_1) = \frac{1}{\sqrt{k}} \max_{1 \leq l \leq n} \left| \sum_{i=1}^{l} I(X_i > X_{(k)}) - \frac{l}{n} \sum_{i=1}^{n} I(X_i > X_{(k)}) \right|, \tag{2.13}$$

which measures the discrepancy between the observed number of the excesses over the high threshold $X_{(k)}$ and the expected number of excesses in each partial time range. Therefore, we reject $\mathcal{H}_0$ if $T_n$ is large. Our analysis shows that $T_n$ has the same limiting distribution as $T_n^\circ$. Based on this, we obtain the result as follows:

**Theorem 1.** *Recall that $\varphi_1(x) = I(x > 0)$ and $\varphi_2(x) = x_+$. Then if the conditions* **(A1)**–**(A3)** *hold, we have that under $\mathcal{H}_0$,*

$$\frac{1}{\sqrt{1+\omega}} T_n(\varphi_1) \Rightarrow \sup_{0 \leq t \leq 1} |B^\circ(t)|, \tag{2.14}$$

*where $B^\circ$ stands for a Brownian bridge.*
  *In addition, if*

$$(k \vee r_n^3) e^{-\varepsilon \sqrt{k}/r_n} = o(1) \qquad \text{for every } \varepsilon > 0, \tag{2.15}$$

*we have that under $\mathcal{H}_0$,*

$$\frac{\alpha}{\sqrt{2+\chi}} T_n(\varphi_2) \Rightarrow \sup_{0 \leq t \leq 1} |B^\circ(t)|. \tag{2.16}$$

**Corollary 1.** *Suppose that under $\mathcal{H}_0$, $\{X_i\}$ is an i.i.d. sequence. Then, under* **(A2)**,

$$T_n(\varphi_1) \Rightarrow \sup_{0 \leq t \leq 1} |B^\circ(t)|$$

*and*

$$\frac{\alpha}{\sqrt{2}} T_n(\varphi_2) \Rightarrow \sup_{0 \leq t \leq 1} |B^\circ(t)|.$$

## 2.2. Cusum test in AR processes

In this subsection, we study the change point test for AR processes. As mentioned earlier in the Introduction, we consider the cusum test based on residuals rather than observa-tions themselves. Let $\{X_i\}$ be an AR($p$) process satisfying the equation:

$$X_i = \sum_{j=1}^{p} \phi_j X_{i-j} + \xi_i,$$

where the characteristic polynomial $\phi(z) := 1 - \phi_1 z - \cdots - \phi_p z^p$ has no zeros inside the unit circle in the complex plane, and $\xi_i$ are error terms. Suppose that one wishes to test

$\mathcal{H}_0^*: \xi_i$ are i.i.d. and the tail index of $\xi_i$ remains the same as for $i = 1, \ldots, n$ vs.

$\mathcal{H}_1^*:$ not $\mathcal{H}_0^*$.

Under $\mathcal{H}_0^*$, we assume that $\{X_i\}$ has a common distribution $F$ and $Z_i := |\xi_i|$ has the distribution $G$. Further, we assume that

**(B1)** $\bar{G} := 1 - G$ is regularly varying at $\infty$, namely,

$$\bar{G}(x) = x^{-\alpha} l^*(x) \qquad \text{for some } \alpha > 0, \tag{2.17}$$

where $l^*$ is slowly varying at $\infty$.

**(B2)** There exist $\kappa(x) = K_\kappa \int_1^x t^{\gamma-1} dt$ ($K_\kappa$ is finite) and a positive measurable function $g \in RV_\gamma$, $\gamma \leq 0$, such that for all $\lambda > 0$,

$$\lim_{x \to \infty} \frac{l^*(\lambda x)/l^*(x) - 1}{g(x)} = \kappa(\lambda).$$

**(B3)** $\sqrt{k} g(b^*(n/k)) \to 0$ as $n \to \infty$, where $b^*(x) = \inf\{y : G(y) \geq 1 - x^{-1}\}$.

For performing a test, we obtain the residuals $\hat{\xi}_i = X_i - \sum_{j=1}^{p} \hat{\phi} X_{i-j}$, where $\hat{\phi} := \hat{\phi}_n$ is an estimator of $\phi = (\phi_1, \ldots, \phi_p)'$ that satisfies the following condition:

**(B4)** There exists a sequence of positive real numbers $\{d(n)\}$ such that

$$d(n) \to \infty \quad \text{and} \quad \frac{\sqrt{k} b^*(n/\sqrt{k})}{b^*(n/k)} = \mathrm{o}(d(n)) \tag{2.18}$$

and

$$d(n)(\hat{\phi} - \phi) = \mathrm{O}_P(1) \qquad \text{as } n \to \infty. \tag{2.19}$$

Typical examples of such $\hat{\phi}$ are the Yule–Walker, linear programming (cf. Feigin and Resnick (1994)), and least gamma deviation (cf. Davis et al. (1992)) estimators.

Then, based on Hill's estimator:

$$\frac{1}{k} \sum_{i=1}^{n} (\log \hat{Z}_i - \log \hat{Z}_{(k+1)})_+, \tag{2.20}$$

where $\hat{Z}_i = |\hat{\xi}_i|$ and $\hat{Z}_{(j)}$ is the $j$th largest r.v. in $\hat{Z}_1, \ldots, \hat{Z}_n$, we employ the cusum test

$$T_n^*(\varphi) := \frac{1}{\sqrt{k}} \max_{1 \leq l \leq n} \left| \sum_{i=1}^{l} \varphi(\log \hat{Z}_i - \log \hat{Z}_{(k)}) - \frac{l}{n} \sum_{i=1}^{n} \varphi(\log \hat{Z}_i - \log \hat{Z}_{(k)}) \right|. \tag{2.21}$$

The following is the main result of this subsection.

**Theorem 2.** *Suppose that conditions* **(B1)**–**(B4)** *hold. Then, under $\mathcal{H}_0^*$, we have*

$$T_n^*(\varphi_1) \Rightarrow \sup_{0 \leq t \leq 1} |B^\circ(t)|. \tag{2.22}$$

*In addition, if*

$$\frac{k^{1/\alpha + 1/2}}{d(n)} = \mathrm{o}(n^{-\nu}) \quad \text{for some } \nu > 0, \tag{2.23}$$

*we have that under $\mathcal{H}_0^*$,*

$$\frac{\alpha}{\sqrt{2}} T_n^*(\varphi_2) \Rightarrow \sup_{0 \leq t \leq 1} |B^\circ(t)|. \tag{2.24}$$

## 2.3. Estimation of the change point under a single abrupt change

In this subsection, we consider the estimating procedure of the location of a change when an abrupt change occurs in the observed time range. We assume that the change point is located at $[n\tau]$ $(0 < \tau < 1)$, and denote the observations by a double array of r.v.'s $X_{n,i}$, $i = 1, \ldots, n$. The common marginal distribution function in time interval $I_{\mathrm{pre}} = \{1, 2, \ldots, [n\tau]\}$ is denoted by $F_{\mathrm{pre}}$ and that in time interval $I_{\mathrm{post}} = \{[n\tau]+1, \ldots, n\}$ is denoted by $F_{\mathrm{post}}$. It is assumed that $\bar{F}_{\mathrm{pre}} \in RV_{-\alpha_{\mathrm{pre}}}$ and $\bar{F}_{\mathrm{post}} \in RV_{-\alpha_{\mathrm{post}}}$.

We set $b_{\mathrm{pre}}(x) = \inf\{y : F_{\mathrm{pre}}(y) \geq 1 - x^{-1}\}$, $b_{\mathrm{post}}(x) = \inf\{y : F_{\mathrm{post}}(y) \geq 1 - x^{-1}\}$, $\mathcal{F}_l^m = \sigma(X_{n,l}, \ldots, X_{n,m})$ and

$$\beta_n(l) = \sup_{m \in \mathbb{N}} \sup\{|P(A \cap B) - P(A)P(B)| : A \in \mathcal{F}_{m+l+1}^n, B \in \mathcal{F}_1^m\}.$$

Here, we focus on the case that $F_{\mathrm{post}}$ has a heavier tail than $F_{\mathrm{pre}}$. The following is the main result of this subsection.

**Theorem 3.** *Suppose that both the $\{X_i = X_{n,i} : i \in I_{\mathrm{pre}}\}$ and $\{X_i = X_{n,i} : i \in I_{\mathrm{post}}\}$ are row-wise stationary. Further, suppose that both the $\bar{F}_{\mathrm{pre}}$ and $\bar{F}_{\mathrm{post}}$ are regularly varying at $\infty$ and there exists a sequence of positive integers $\{r_n\}$ such that*

$$m_n \beta_n(r_n) = \mathrm{o}(1) \quad \text{and} \quad r_n = \mathrm{o}(k), \tag{2.25}$$

where $m_n = [n/r_n]$. Then if there exists $c > 1$ such that

$$\liminf_{x \to \infty} \frac{b_{\text{post}}(x)}{b_{\text{pre}}(x)} > c, \tag{2.26}$$

we have $T_n(\varphi_1) \xrightarrow{P} \infty$; further, if there exists $d \in (1, \infty]$ such that

$$\lim_{x \to \infty} \frac{b_{\text{post}}(x)}{b_{\text{pre}}(x)} = d, \tag{2.27}$$

then $T_n(\varphi_2) \xrightarrow{P} \infty$ and $\hat{\tau}_n := \hat{l}_n/n$ with

$$\hat{l}_n := \arg\max_{1 \leq l \leq n} \left| \sum_{i=1}^{l} \varphi(\log X_i - \log X_{(k)}) - \frac{l}{n} \sum_{i=1}^{n} \varphi(\log X_i - \log X_{(k)}) \right|$$

is a consistent estimator of $\tau$, namely,

$$\hat{\tau}_n \xrightarrow{P} \tau. \tag{2.28}$$

## 3. Simulation study

In this section, we evaluate the performance of the proposed tests through a simulation study. Here, we employ the decision rule: At the nominal level 0.05 we reject $\mathcal{H}_0$ if the scaled $T_n(\varphi)$ is greater than 1.35, the scaling constant of which depends on $\varphi$ and the dependency of data. Given any significance levels, the critical values can be obtained from a Monte Carlo simulation (cf. Lee *et al.* (2003)). In what follows, we briefly explain how to obtain those. We generate the random numbers $\varepsilon_1, \ldots, \varepsilon_N$, $N = 10\,000$, following the standard normal distribution, and calculate

$$L := \frac{1}{\sqrt{N}} \max_{1 \leq l \leq N} \left| \sum_{i=1}^{l} \varepsilon_i - \frac{l}{N} \sum_{i=1}^{N} \varepsilon_i \right|.$$

Then we determine the critical value at the nominal level 0.05 as the 0.95-quantile from such 10 000 $L$'s. Table 1 presents the critical values for the nominal levels 0.01, 0.05 and 0.1.

In our simulation study, we consider the following distributions:

**Table 1.** Critical values

| Nominal level  | 0.9  | 0.95 | 0.99 |
|---|---|---|---|
| Critical value | 1.22 | 1.35 | 1.60 |

- Burr distribution for the i.i.d. sample:

$$\bar{F}(x) = \left(\frac{\beta}{\beta + x^{-\gamma}}\right)^{\lambda} \qquad (\lambda > 0, \beta > 0, \gamma < 0).$$

Its tail index is the reciprocal of $\alpha = -\gamma\lambda$ and the second-order regularly varying exponent is $\gamma$. This is mainly used to investigate the effect of the second-order regularly varying exponent $\gamma$ on the tests.

- $t$ distribution with the $\nu$ degrees of freedom for AR and moving average (MA) models. Its tail index is $\nu^{-1}$, that is, $\alpha = \nu$ and the second-order regularly varying exponent is $\gamma = -2/\nu$. This is used to generate the innovations of the moving average and autoregressive processes.

We first consider the i.i.d. case with the Burr distribution. In this case, we reject $\mathcal{H}_0$ if $T_n(\varphi_1) \geq 1.35$, and also if $\frac{\hat{\alpha}}{\sqrt{2}} T_n(\varphi_2) \geq 1.35$ with

$$\hat{\alpha} = \left\{\frac{1}{k}\sum_{i=1}^{n}(\log X_i - \log X_{(k+1)})_+\right\}^{-1}.$$

Tables 2 and 3 show that the empirical sizes of $T_n(\varphi_1)$ are closer to the nominal level than $T_n(\varphi_2)$, which means the former is more stable than the latter. This phenomenon can be seen in all other cases considered in this simulation. From Table 2, it can be seen that $\alpha$ and $\gamma$ do not affect the performance of $T_n(\varphi_1)$ much. On the other hand, from Table 3, it can be seen that $\gamma$ affects the performance of $T_n(\varphi_2)$ to certain degree; the empirical size becomes smaller as $\gamma$ gets close to 0 while this phenomenon is not elaborate with the change in $\alpha$. In fact, $T_n(\varphi_1)$ depends on the rank and observed time of observations but not on their magnitude.

Next, we deal with the MA(1) process $X_i = \xi_i + \theta\xi_{i-1}$, where $\xi_i$ are i.i.d. innovations following a $t$ distribution with the $\nu$ degrees of freedom. To perform a test, we estimate $\chi$ and $\omega$ (cf. Theorem 1) by the estimators

$$\hat{\chi} = \frac{2\hat{\alpha}}{k}\sum_{i=1}^{n-1}(\log X_i - \log X_{(k)})_+(\log X_{i+1} - \log X_{(k)})_+$$

and

$$\hat{\omega} = \frac{2}{k}\sum_{i=1}^{n-1} I(X_i > X_{(k)}, X_{i+1} > X_{(k)})$$

(cf. Hsing (1991)). In this case, we reject $\mathcal{H}_0$ if $\frac{1}{\sqrt{1+\hat{\omega}}}T_n(\varphi_1) \geq 1.35$, and also if $\frac{\hat{\alpha}}{\sqrt{2+\hat{\chi}}}T_n(\varphi_2) \geq 1.35$. Tables 4 and 5 exhibit the empirical sizes of $T_n(\varphi_1)$ and $T_n(\varphi_2)$, respectively. Although some size distortions exist, particularly when $n$ is 1000, this size distortion effect seems to be soothed as $n$ increases.

**Table 2.** Empirical sizes of $T_n(\varphi_1)$ in the i.i.d. case

| | | | $k$ | | | | | | | | | |
|---|---|---|---|---|---|---|---|---|---|---|---|---|
| $n$ | $\alpha$ | $\gamma$ | 10 | 20 | 30 | 40 | 50 | 60 | 70 | 80 | 90 | 100 |
| 1000 | 2 | $-2.0$ | 0.035 | 0.040 | 0.035 | 0.037 | 0.035 | 0.037 | 0.033 | 0.032 | 0.034 | 0.033 |
| | 2 | $-0.5$ | 0.030 | 0.041 | 0.039 | 0.035 | 0.038 | 0.035 | 0.037 | 0.036 | 0.032 | 0.031 |
| | 1 | $-2.0$ | 0.035 | 0.040 | 0.037 | 0.041 | 0.036 | 0.033 | 0.033 | 0.033 | 0.031 | 0.028 |
| | 1 | $-0.5$ | 0.030 | 0.036 | 0.038 | 0.038 | 0.037 | 0.036 | 0.034 | 0.029 | 0.031 | 0.032 |
| | | | $k$ | | | | | | | | | |
| $n$ | $\alpha$ | $\gamma$ | 25 | 50 | 75 | 100 | 125 | 150 | 175 | 200 | 225 | 250 |
| 3000 | 2 | $-2.0$ | 0.033 | 0.048 | 0.043 | 0.044 | 0.038 | 0.033 | 0.024 | 0.043 | 0.037 | 0.039 |
| | 2 | $-0.5$ | 0.031 | 0.047 | 0.038 | 0.040 | 0.035 | 0.033 | 0.033 | 0.035 | 0.044 | 0.038 |
| | 1 | $-2.0$ | 0.045 | 0.047 | 0.031 | 0.046 | 0.035 | 0.035 | 0.042 | 0.043 | 0.029 | 0.037 |
| | 1 | $-0.5$ | 0.033 | 0.041 | 0.040 | 0.034 | 0.056 | 0.038 | 0.048 | 0.034 | 0.046 | 0.039 |

Now, we turn our attention to the AR(1) process $X_i = \phi X_{i-1} + \xi_i$, where $\xi_i$ are identical to those in the previous case. In order to perform the tests based on residuals, we employ $\phi = 0.5$ and $0.9$ and obtain the residuals $\hat{\xi}_i = X_i - \hat{\phi} X_{i-1}$ by using the least squares estimator. According to Theorem 2, the decision rules are the same as those in the i.i.d. case, but $T_n(\varphi)$ and $\mathcal{H}_0$ are replaced by $T_n^*(\varphi)$ and $\mathcal{H}_0^*$, respectively. Tables 6 and 7 show that the empirical sizes are as good as those in the i.i.d. case regardless of the value of $\phi$'s.

So far, we have investigated the stability of the tests. In what follows, we examine the power of $T_n(\varphi_1)$ for the MA(1) and AR(1) processes and the associated MSE of $\hat{\tau}_n$ for

**Table 3.** Empirical sizes of the $T_n(\varphi_2)$ in the i.i.d. case

| | | | $k$ | | | | | | | | | |
|---|---|---|---|---|---|---|---|---|---|---|---|---|
| $n$ | $\alpha$ | $\gamma$ | 10 | 20 | 30 | 40 | 50 | 60 | 70 | 80 | 90 | 100 |
| 1000 | 2 | $-2.0$ | 0.011 | 0.021 | 0.029 | 0.028 | 0.029 | 0.032 | 0.029 | 0.028 | 0.027 | 0.032 |
| | 2 | $-0.5$ | 0.009 | 0.017 | 0.022 | 0.021 | 0.023 | 0.021 | 0.019 | 0.021 | 0.018 | 0.019 |
| | 1 | $-2.0$ | 0.012 | 0.023 | 0.029 | 0.029 | 0.030 | 0.031 | 0.031 | 0.032 | 0.035 | 0.031 |
| | 1 | $-0.5$ | 0.009 | 0.019 | 0.025 | 0.027 | 0.023 | 0.025 | 0.024 | 0.026 | 0.024 | 0.023 |
| | | | $k$ | | | | | | | | | |
| $n$ | $\alpha$ | $\gamma$ | 25 | 50 | 75 | 100 | 125 | 150 | 175 | 200 | 225 | 250 |
| 3000 | 2 | $-2.0$ | 0.023 | 0.033 | 0.033 | 0.041 | 0.037 | 0.032 | 0.029 | 0.039 | 0.020 | 0.039 |
| | 2 | $-0.5$ | 0.025 | 0.027 | 0.031 | 0.034 | 0.037 | 0.030 | 0.031 | 0.021 | 0.027 | 0.030 |
| | 1 | $-2.0$ | 0.031 | 0.026 | 0.032 | 0.039 | 0.044 | 0.040 | 0.043 | 0.034 | 0.035 | 0.029 |
| | 1 | $-0.5$ | 0.029 | 0.031 | 0.041 | 0.038 | 0.035 | 0.032 | 0.038 | 0.029 | 0.028 | 0.035 |

**Table 4.** Empirical sizes of the $T_n(\varphi_1)$ for the MA(1) model

| $n$ | $\alpha$ | $\gamma$ | $k$=10 | 20 | 30 | 40 | 50 | 60 | 70 | 80 | 90 | 100 |
|---|---|---|---|---|---|---|---|---|---|---|---|---|
| 1000 | 2 | 0.1 | 0.030 | 0.031 | 0.028 | 0.026 | 0.030 | 0.016 | 0.021 | 0.019 | 0.021 | 0.015 |
|  | 2 | 0.5 | 0.031 | 0.030 | 0.032 | 0.030 | 0.025 | 0.028 | 0.023 | 0.021 | 0.024 | 0.019 |
|  | 2 | 1.0 | 0.031 | 0.026 | 0.026 | 0.030 | 0.026 | 0.025 | 0.029 | 0.024 | 0.022 | 0.022 |

| $n$ | $\alpha$ | $\theta$ | $k$=25 | 50 | 75 | 100 | 125 | 150 | 175 | 200 | 225 | 250 |
|---|---|---|---|---|---|---|---|---|---|---|---|---|
| 3000 | 2 | 0.1 | 0.047 | 0.038 | 0.032 | 0.027 | 0.024 | 0.039 | 0.025 | 0.026 | 0.013 | 0.019 |
|  | 2 | 0.5 | 0.034 | 0.039 | 0.034 | 0.036 | 0.039 | 0.035 | 0.020 | 0.042 | 0.030 | 0.022 |
|  | 2 | 1.0 | 0.025 | 0.043 | 0.046 | 0.044 | 0.034 | 0.034 | 0.032 | 0.035 | 0.024 | 0.032 |

**Table 5.** Empirical sizes of the $T_n(\varphi_2)$ for the MA(1) model

| $n$ | $\alpha$ | $\gamma$ | $k$=10 | 20 | 30 | 40 | 50 | 60 | 70 | 80 | 90 | 100 |
|---|---|---|---|---|---|---|---|---|---|---|---|---|
| 1000 | 2 | 0.1 | 0.004 | 0.022 | 0.018 | 0.023 | 0.019 | 0.017 | 0.022 | 0.018 | 0.011 | 0.016 |
|  | 2 | 0.5 | 0.004 | 0.010 | 0.015 | 0.017 | 0.017 | 0.019 | 0.020 | 0.022 | 0.015 | 0.019 |
|  | 2 | 1.0 | 0.002 | 0.012 | 0.016 | 0.021 | 0.017 | 0.018 | 0.019 | 0.022 | 0.013 | 0.019 |

| $n$ | $\alpha$ | $\theta$ | $k$=25 | 50 | 75 | 100 | 125 | 150 | 175 | 200 | 225 | 250 |
|---|---|---|---|---|---|---|---|---|---|---|---|---|
| 3000 | 2 | 0.1 | 0.030 | 0.022 | 0.029 | 0.035 | 0.028 | 0.033 | 0.029 | 0.021 | 0.029 | 0.031 |
|  | 2 | 0.5 | 0.007 | 0.015 | 0.029 | 0.025 | 0.022 | 0.034 | 0.026 | 0.032 | 0.028 | 0.029 |
|  | 2 | 1.0 | 0.020 | 0.030 | 0.023 | 0.030 | 0.033 | 0.038 | 0.033 | 0.028 | 0.027 | 0.016 |

**Table 6.** Empirical sizes of $T_n^*(\varphi_1)$ for the AR(1) model

| $n$ | $\alpha$ | $\phi$ | $k$=10 | 20 | 30 | 40 | 50 | 60 | 70 | 80 | 90 | 100 |
|---|---|---|---|---|---|---|---|---|---|---|---|---|
| 1000 | 2 | 0.5 | 0.034 | 0.037 | 0.037 | 0.039 | 0.035 | 0.035 | 0.033 | 0.035 | 0.030 | 0.030 |
|  | 2 | 0.9 | 0.032 | 0.035 | 0.037 | 0.036 | 0.036 | 0.039 | 0.036 | 0.033 | 0.036 | 0.033 |

| $n$ | $\alpha$ | $\phi$ | $k$=25 | 50 | 75 | 100 | 125 | 150 | 175 | 200 | 225 | 250 |
|---|---|---|---|---|---|---|---|---|---|---|---|---|
| 3000 | 2 | 0.5 | 0.045 | 0.046 | 0.037 | 0.049 | 0.044 | 0.043 | 0.033 | 0.046 | 0.034 | 0.043 |
|  | 2 | 0.9 | 0.031 | 0.049 | 0.046 | 0.039 | 0.039 | 0.040 | 0.037 | 0.032 | 0.041 | 0.039 |

**Table 7.** Empirical sizes of $T_n^*(\varphi_2)$ for the AR(1) model

| | | | $k$ | | | | | | | | | |
|---|---|---|---|---|---|---|---|---|---|---|---|---|
| $n$ | $\alpha$ | $\phi$ | 10 | 20 | 30 | 40 | 50 | 60 | 70 | 80 | 90 | 100 |
| 1000 | 2 | 0.5 | 0.012 | 0.021 | 0.026 | 0.024 | 0.028 | 0.029 | 0.024 | 0.026 | 0.028 | 0.022 |
| | 2 | 0.9 | 0.010 | 0.022 | 0.026 | 0.023 | 0.024 | 0.026 | 0.028 | 0.029 | 0.024 | 0.021 |
| | | | $k$ | | | | | | | | | |
| $n$ | $\alpha$ | $\phi$ | 25 | 50 | 75 | 100 | 125 | 150 | 175 | 200 | 225 | 250 |
| 3000 | 2 | 0.5 | 0.024 | 0.034 | 0.034 | 0.033 | 0.024 | 0.033 | 0.027 | 0.039 | 0.029 | 0.032 |
| | 2 | 0.9 | 0.031 | 0.032 | 0.037 | 0.039 | 0.040 | 0.032 | 0.030 | 0.027 | 0.027 | 0.019 |

**Table 8.** Empirical powers of $T_n(\varphi_1)$ for the MA(1) model

| | | $k$ | | | | | | | | | |
|---|---|---|---|---|---|---|---|---|---|---|---|
| $n$ | $\tau$ | 10 | 20 | 30 | 40 | 50 | 60 | 70 | 80 | 90 | 100 |
| 1000 | 0.25 | 0.06 | 0.15 | 0.29 | 0.46 | 0.71 | 0.81 | 0.85 | 0.84 | 0.84 | 0.84 |
| | 0.50 | 0.44 | 0.99 | 0.99 | 0.99 | 0.99 | 0.99 | 0.99 | 0.99 | 0.98 | 0.99 |
| | 0.75 | 0.93 | 0.97 | 0.97 | 0.97 | 0.95 | 0.94 | 0.93 | 0.89 | 0.85 | 0.82 |
| | | $k$ | | | | | | | | | |
| $n$ | $\tau$ | 25 | 50 | 75 | 100 | 125 | 150 | 175 | 200 | 225 | 250 |
| 3000 | 0.25 | 0.24 | 0.94 | 0.99 | 1.00 | 1.00 | 1.00 | 1.00 | 1.00 | 1.00 | 1.00 |
| | 0.50 | 1.00 | 1.00 | 1.00 | 1.00 | 1.00 | 1.00 | 1.00 | 1.00 | 1.00 | 1.00 |
| | 0.75 | 1.00 | 1.00 | 1.00 | 1.00 | 1.00 | 1.00 | 1.00 | 1.00 | 1.00 | 1.00 |

the MA(1) process; $T_n(\varphi_2)$ is not considered here since its performance has a pattern similar to that of $T_n(\varphi_1)$ and is not so good compared to $T_n(\varphi_1)$. To task this, we take account of the MA(1) process $X_i = \xi_i + 0.5\xi_{i-1}$ and the AR(1) process $X_i = 0.5X_{i-1} + \xi_i$. Under the alternative hypothesis, it is assumed that $\xi_i \sim t(3)$ for $i \leq [n\tau]$ and $\xi_i \sim t(1)$ for $i > [n\tau]$ (i.e., $\nu$ changes from 3 to 1) with $\tau = 0.25, 0.5, 0.75$. Tables 8 and 9 exhibit the empirical powers of $T_n(\varphi_1)$ and $T_n^*(\varphi_1)$, respectively, and show that these tests produce reasonably good powers. As might be anticipated, the tests have a tendency to produce the best powers when $\tau = 0.5$. Table 10 also shows that the same phenomenon is true for the MSE of $\hat{\tau}_n$.

From Tables 2–9, it can be seen that the performance of the tests does not much depend on the tail sample fraction $k$. Our simulation result indicates that it is not an easy task to build up a rule to choose an optimal tail sample fraction unlike the case of Hill's estimator since we cannot see an obvious trend from our results. Our findings only enable us to recommend the use of $k$ within a reasonable range, say, $k = 50$–$100$ when $n = 1000$ and $k = 100\text{--}200$ when $n = 3000$. Due to its importance, we leave the issue of finding an optimal $k$ as a task for future study.

**Table 9.** Empirical powers of $T_n^*(\varphi_1)$ for the AR(1) model

| | | k | | | | | | | | | |
|---|---|---|---|---|---|---|---|---|---|---|---|
| $n$ | $\tau$ | 10 | 20 | 30 | 40 | 50 | 60 | 70 | 80 | 90 | 100 |
| 1000 | 0.25 | 0.11 | 0.34 | 0.75 | 0.86 | 0.90 | 0.89 | 0.90 | 0.89 | 0.86 | 0.83 |
| | 0.50 | 0.94 | 0.99 | 0.99 | 0.99 | 0.99 | 0.98 | 0.99 | 0.98 | 0.98 | 0.97 |
| | 0.75 | 0.94 | 0.97 | 0.98 | 0.96 | 0.95 | 0.92 | 0.91 | 0.87 | 0.83 | 0.78 |
| | | k | | | | | | | | | |
| $n$ | $\tau$ | 25 | 50 | 75 | 100 | 125 | 150 | 175 | 200 | 225 | 250 |
| 3000 | 0.25 | 0.65 | 1.00 | 1.00 | 1.00 | 1.00 | 1.00 | 1.00 | 1.00 | 1.00 | 1.00 |
| | 0.50 | 1.00 | 1.00 | 1.00 | 1.00 | 1.00 | 1.00 | 1.00 | 1.00 | 1.00 | 1.00 |
| | 0.75 | 1.00 | 1.00 | 1.00 | 1.00 | 1.00 | 1.00 | 1.00 | 1.00 | 1.00 | 1.00 |

**Table 10.** MSE of $\hat{\tau}_n$ of $T_n(\varphi_1)$ for the MA(1) model

| | | k | | | | | | | | | |
|---|---|---|---|---|---|---|---|---|---|---|---|
| $n$ | $\tau$ | 10 | 20 | 30 | 40 | 50 | 60 | 70 | 80 | 90 | 100 |
| 1000 | 0.25 | 0.090 | 0.056 | 0.036 | 0.034 | 0.024 | 0.024 | 0.019 | 0.017 | 0.018 | 0.016 |
| | 0.50 | 0.019 | 0.007 | 0.004 | 0.003 | 0.003 | 0.003 | 0.002 | 0.002 | 0.002 | 0.002 |
| | 0.75 | 0.003 | 0.001 | 0.002 | 0.002 | 0.003 | 0.003 | 0.003 | 0.005 | 0.007 | 0.007 |
| | | k | | | | | | | | | |
| $n$ | $\tau$ | 25 | 50 | 75 | 100 | 125 | 150 | 175 | 200 | 225 | 251 |
| 3000 | 0.25 | 0.047 | 0.025 | 0.016 | 0.010 | 0.007 | 0.006 | 0.005 | 0.005 | 0.004 | 0.004 |
| | 0.50 | 0.005 | 0.002 | 0.001 | 0.000 | 0.000 | 0.000 | 0.000 | 0.000 | 0.000 | 0.000 |
| | 0.75 | 0.000 | 0.000 | 0.000 | 0.000 | 0.000 | 0.000 | 0.000 | 0.000 | 0.000 | 0.000 |

## 4. Proofs

**Lemma 1.** *Under Condition* **(A2)**, *we have that for every* $0 < v < 1 < u < \infty$,

$$\lim_{x \to \infty} \sup_{\lambda \in [v,u] \setminus \{1\}} \left| \frac{l(\lambda x)/l(x) - 1}{\kappa(\lambda) g(x)} - 1 \right| = 0. \tag{4.1}$$

**Proof.** According to Theorem 2.3 of Hsing (1991), we have

$$\lim_{x\to\infty}\sup_{\lambda\in(1,u\vee v^{-1}]}\left|\frac{l(\lambda x)/l(x)-1}{\kappa(\lambda)g(x)}-1\right|=0. \tag{4.2}$$

Further, by (4.2) and the fact that $\lambda^{-\gamma}\kappa(\lambda)=-\kappa(\lambda^{-1})$,

$$\lim_{x\to\infty}\sup_{\lambda\in[v,1)}\left|\frac{l(\lambda x)}{l(x)}-1\right|=0,$$

and

$$\lim_{x\to\infty}\sup_{\lambda\in[v,1)}\left|\frac{g(\lambda x)}{\lambda^\gamma g(x)}-1\right|=0$$

(cf. Theorem 1.2.1 in Bingham *et al.* (1987)), we have

$$\lim_{x\to\infty}\sup_{\lambda\in[v,1)}\left|\frac{l(\lambda x)/l(x)-1}{g(x)\kappa(\lambda)}-1\right|=\lim_{x\to\infty}\sup_{\lambda\in[v,1)}\left|\frac{l(\lambda x)}{l(x)}\frac{g(\lambda x)}{\lambda^\gamma g(x)}\frac{l(x)/l(\lambda x)-1}{g(\lambda x)\kappa(\lambda^{-1})}-1\right|=0.$$

Combining this and (4.2), we assert (4.1). □

**Lemma 2.** *Under Condition* **(A2)**,

$$\bar{F}(e^{\zeta/\sqrt{k}}b(n/k))=\frac{k}{n}\left(1-\frac{\alpha\zeta}{\sqrt{k}}+o\left(\frac{1}{\sqrt{k}}\right)\right) \tag{4.3}$$

*uniformly on every compact $\zeta$-set in $\mathbb{R}$. Further, for any $r>0$,*

$$\begin{aligned}\mathrm{E}&\left|\left(\log X_1-\log b(n/k)+\frac{\zeta_1}{\sqrt{k}}\right)_+-\left(\log X_1-\log b(n/k)+\frac{\zeta_2}{\sqrt{k}}\right)_+\right|^r\\&\leq\left|\frac{\zeta_1-\zeta_2}{\sqrt{k}}\right|^r\frac{k}{n}\left(1+\mathrm{O}\left(\frac{1}{\sqrt{k}}\right)\right)\end{aligned} \tag{4.4}$$

*uniformly on every compact $(\zeta_1,\zeta_2)$-set in $\mathbb{R}^2$.*

**Proof.** We first verify (4.3). According to the arguments in the proof of Theorem 2.4 of Hsing (1991), we can express

$$\bar{F}(b(n/k))=\frac{k}{n}\left(1+o\left(\frac{1}{\sqrt{k}}\right)\right).$$

Hence, we can express that for $\zeta\in[-K,K]$, $K>0$,

$$\bar{F}(e^{\zeta/\sqrt{k}}b(n/k))=\bar{F}(b(n/k))\frac{\bar{F}(e^{\zeta/\sqrt{k}}b(n/k))}{\bar{F}(b(n/k))}$$

$$= \frac{k}{n}\left(1+\mathrm{o}\left(\frac{1}{\sqrt{k}}\right)\right)\mathrm{e}^{-\alpha\zeta/\sqrt{k}}\frac{l(\mathrm{e}^{\zeta/\sqrt{k}}b(n/k))}{l(b(n/k))}$$

$$= \frac{k}{n}\left(1+\mathrm{o}\left(\frac{1}{\sqrt{k}}\right)\right)\left(1-\frac{\alpha\zeta}{\sqrt{k}}+\Delta_{n,1}(\zeta)\right)\frac{l(\mathrm{e}^{\zeta/\sqrt{k}}b(n/k))}{l(b(n/k))},$$

where

$$\lim_{n\to\infty}\sup_{\zeta\in[-K,K]}|\sqrt{k}\Delta_{n,1}(\zeta)|=0.$$

Due to Lemma 1, we have

$$\frac{l(\mathrm{e}^{\zeta/\sqrt{k}}b(n/k))}{l(b(n/k))}=1+\kappa(\mathrm{e}^{\zeta/\sqrt{k}})g(b(n/k))(1+\Delta_{n,2}(\zeta)),$$

where

$$\lim_{n\to\infty}\sup_{\zeta\in[-K,K]}\max\{|\Delta_{n,2}(\zeta)|,|\kappa(\mathrm{e}^{\zeta/\sqrt{k}})|\}=0.$$

Henceforth, since $\lim_{n\to\infty}\sqrt{k}g(b(n/k))<\infty$, we have

$$\lim_{n\to\infty}\sup_{\zeta\in[-K,K]}\sqrt{k}\left|\frac{l(\mathrm{e}^{\zeta/\sqrt{k}}b(n/k))}{l(b(n/k))}-1\right|=0,$$

which asserts (4.3). Now that (4.4) can be proven by using (4.3), the inequality: $|x_+ - y_+|^r \leq |x-y|^r I(\max\{x,y\}>0)$, and the fact that $P(\log X_1 - \log b(n/k) + \frac{\zeta}{\sqrt{k}} > 0) = \bar{F}(\mathrm{e}^{-\zeta/\sqrt{k}}b(n/k))$, the lemma is established. □

For $\zeta \in \mathbb{R}$, we set

$$A_i(\zeta):=A_{ni}(\zeta,\varphi):=\varphi\left(\log X_i - \log b(n/k) + \frac{\zeta}{\sqrt{k}}\right)$$

and

$$M_n(t,\zeta):=M_n(t,\zeta,\varphi):=\frac{1}{\sqrt{k}}\sum_{i=1}^{[nt]}\{A_i(\zeta)-\mathrm{E}A_1(\zeta)\}.$$

Note that $M_n(t)$ in Section 2.1 is identical to $M_n(t,0)$.

**Lemma 3.** *Suppose that conditions* **(A1)**–**(A3)** *hold. Then, for every $\zeta \in \mathbb{R}$ and $t \in [0,1]$,*

$$M_n(t,\zeta)-M_n(t,0)=\mathrm{o}_P(1). \tag{4.5}$$

*Further,*

$$\sup_{0\leq t\leq 1} |M_n(t,\zeta) - M_n(t,0)| = o_P(1). \qquad (4.6)$$

**Proof.** In order to verify (4.5), it suffices to show that for any real $\zeta$ and $t \in (0,1)$,

$$\frac{1}{\sqrt{k}} \sum_{i=1}^{[nt]} \{A_i(0) - A_i(\zeta) - E(A_1(0) - A_1(\zeta))\} = o_P(1). \qquad (4.7)$$

Let $\{r_n\}$ be a sequence that satisfies Condition **(A1)**, and let $m_n$ be the integer part of $[nt]/r_n$. To show (4.7), we express the left-hand side of (4.7) as $\sum_{i=1}^{m_n} S_{ni} + R_n$, where

$$S_{ni} = \frac{1}{\sqrt{k}} \sum_{j=(i-1)r_n+1}^{ir_n} \{A_j(0) - A_j(\zeta) - E(A_1(0) - A_1(\zeta))\}$$

and

$$R_n = \frac{1}{\sqrt{k}} \sum_{j=m_n r_n+1}^{[nt]} \{A_j(0) - A_j(\zeta) - E(A_1(0) - A_1(\zeta))\}.$$

We first verify that

$$\sum_{i \in O_n} S_{ni} = o_P(1) \qquad \text{as } n \to \infty, \qquad (4.8)$$

where $O_n$ denotes the set of the odd numbers in $\{1, 2, \ldots, m_n\}$. Let $\{\tilde{S}_{ni} : i \in O_n\}$ be i.i.d. copies of $S_{n1}$. Note that for each $\varepsilon > 0$,

$$P\left(\left|\sum_{i\in O_n} \tilde{S}_{ni}\right| > \varepsilon\right) \leq \frac{1}{\varepsilon^2} \sum_{i\in O_n} \operatorname{Var} S_{ni} \leq \frac{[nt]r_n}{2k\varepsilon^2} E(A_1(0) - A_1(\zeta))^2 = o(1)$$

owing to (2.8) and the fact that $nE(A_1(0) - A_1(\zeta))^2 = O(\sqrt{k})$ (cf. Lemma 2). This implies

$$\sum_{i\in O_n} \tilde{S}_{ni} = o_P(1) \qquad \text{as } n \to \infty. \qquad (4.9)$$

Further, according to Lemma 2 of Billingsley (1995), page 365 and (2.7), it can be yielded that for each $s \in \mathbb{R}$,

$$\left| E\exp\left\{s\mathbf{i} \sum_{i\in O_{ni}} S_{ni}\right\} - \prod_{i\in O_{ni}} E\exp\{s\mathbf{i} S_{ni}\} \right| \leq 16 m_n \beta(r_n) \to 0 \qquad \text{as } n \to \infty.$$

Thus, from this and (4.9), (4.8) is obtained. In fact, a similar result can be obtained for the summation of $S_{ni}$'s over $E_n := \{1, \ldots, m_n\} \backslash O_n$. Therefore, we have $\sum_{i=1}^{m_n} S_{ni} = o_P(1)$.

Since $R_n$ is also $o_P(1)$ by the fact:

$$ER_n^2 \leq \frac{r_n^2}{k} E(A_1(0) - A_1(\zeta))^2 = o(1),$$

which is due to (2.8) and Lemma 2, (4.5) is asserted.

Now, to verify (4.6) we only need to show that

$$\lim_{\rho \to 0} \limsup_n P\left(\sup_{|t_1-t_2|<\rho} |M_n(t_1,\zeta) - M_n(t_1,0) - \{M_n(t_2,\zeta) - M_n(t_2,0)\}| > \varepsilon\right) = 0 \tag{4.10}$$

owing to (4.5). To task this, note that for $\zeta \in \mathbb{R}$, $\rho > 0$ and $0 \leq t_1 \leq t_2 \leq 1$,

$$\sup_{|t_1-t_2|<\rho} |M_n(t_1,\zeta) - M_n(t_1,0) - \{M_n(t_2,\zeta) - M_n(t_2,0)\}|$$

$$= \sup_{|t_1-t_2|<\rho} \frac{1}{\sqrt{k}} \left| \sum_{i=[nt_1]+1}^{[nt_2]} (A_i(0) - A_i(\zeta) - \mathrm{E}\{A_1(0) - A_1(\zeta)\}) \right|$$

$$\leq \sup_{|t_1-t_2|<\rho} \frac{1}{\sqrt{k}} \sum_{i=[nt_1]+1}^{[nt_2]} |A_i(0) - A_i(\zeta)| + \sup_{|t_1-t_2|\leq\rho} \frac{[nt_2] - [nt_1]}{\sqrt{k}} \mathrm{E}|A_1(0) - A_1(\zeta)|$$

$$= I_n(\rho) + II_n(\rho).$$

It can be seen that

$$I_n(\rho) \leq \max_{0 \leq j \leq [\rho^{-1}]-1} \frac{1}{\sqrt{k}} \sum_{i=[ns_j]+1}^{[ns_{j+2}]} |A_i(0) - A_i(\zeta)|,$$

where $s_j = j\rho \wedge 1$, and thus, for any $\varepsilon > 0$,

$$P(I_n(\rho) \geq \varepsilon) \leq \frac{1}{\rho} P\left(\frac{1}{\sqrt{k}} \sum_{i=1}^{[2\rho n]+1} |A_i(0) - A_i(\zeta)| > \varepsilon\right).$$

According to Lemma 2, we can choose a sufficiently small $\rho_0 > 0$ such that

$$\limsup_{n \to \infty} \frac{[2\rho_0 n]+1}{\sqrt{k}} \mathrm{E}|A_1(0) - A_1(\zeta)| < \varepsilon/2.$$

Therefore, due to (4.5), for any $\rho \leq \rho_0$,

$$\limsup_{n \to \infty} P\left(\frac{1}{\sqrt{k}} \sum_{i=1}^{[2\rho n]+1} |A_i(0) - A_i(\zeta)| > \varepsilon\right)$$

$$\leq \limsup_{n\to\infty} P\left(\frac{1}{\sqrt{k}} \sum_{i=1}^{[2\rho n]+1} \{|A_i(0) - A_i(\zeta)| - \mathrm{E}|A_1(0) - A_1(\zeta)|\} > \frac{\varepsilon}{2}\right) = 0,$$

which yields $\lim_{\rho\to 0} \limsup_n P(I_n(\rho) > \varepsilon) = 0$. Since $II_n(\rho) \leq \rho \mathrm{O}(1)$ as $n \to \infty$, (4.6) is asserted. This completes the proof. □

**Lemma 4.** *Suppose that conditions* **(A1)–(A3)** *hold. Then, for every $\varepsilon > 0$ and $K > 0$,*

$$\lim_{\rho\to 0} \limsup_{n} P\left(\sup_{|\zeta_1-\zeta_2|<\rho} \frac{1}{\sqrt{k}} \sum_{i=1}^{n} |A_i(\zeta_2) - A_i(\zeta_1)| > \varepsilon\right) = 0 \tag{4.11}$$

*and*

$$\lim_{\rho\to 0} \limsup_{n} \sup_{|\zeta_1-\zeta_2|<\rho} \frac{n}{\sqrt{k}} \mathrm{E}|A_1(\zeta_2) - A_1(\zeta_1)| = 0, \tag{4.12}$$

*where both $\zeta_1$ and $\zeta_2$ are numbers in $[-K, K]$.*

**Proof.** Since (4.12) can be directly obtained from Lemma 2, we only prove (4.11). For simplicity, we assume that $2K/\rho$ is an integer. Note that for $\zeta_1$ and $\zeta_2 \in [-K, K]$,

$$\sup_{|\zeta_1-\zeta_2|<\rho} \frac{1}{\sqrt{k}} \sum_{i=1}^{n} |A_i(\zeta_2) - A_i(\zeta_1)| \leq \max_{l} \frac{1}{\sqrt{k}} \sum_{i=1}^{n} \{A_i((l+2)\rho) - A_i(l\rho)\}, \tag{4.13}$$

where the maximum is taken over the integers $l$ such that $[l\rho, (l+2)\rho] \subset [-K, K]$. From (4.5), it can be easily seen that

$$\frac{1}{\sqrt{k}} \sum_{i=1}^{n} \{A_i((l+2)\rho) - A_i(l\rho) - \mathrm{E}(A_1((l+2)\rho) - A_1(l\rho))\} = \mathrm{o}_P(1).$$

By using this, (4.12) and (4.13), we can verify (4.11) similarly to Lemma 3. □

**Lemma 5.** *Suppose that conditions* **(A1)–(A3)** *hold. Then, for every $K > 0$,*

$$\sup_{\zeta \in [-K,K]} \sup_{0 \leq t \leq 1} |M_n(t,\zeta) - M_n(t,0)| = \mathrm{o}_P(1).$$

**Proof.** Let $\varepsilon > 0$. Note that for $\zeta_1$ and $\zeta_2 \in [-K, K]$,

$$\sup_{|\zeta_1-\zeta_2|<\rho} \sup_{0\leq t\leq 1} |M_n(t,\zeta_1) - M_n(t,\zeta_2)|$$

$$= \sup_{|\zeta_1-\zeta_2|<\rho} \sup_{0\leq t\leq 1} \left|\frac{1}{\sqrt{k}} \sum_{i=1}^{[nt]} \{A_i(\zeta_1) - \mathrm{E}A_1(\zeta_1)\} - \frac{1}{\sqrt{k}} \sum_{i=1}^{[nt]} \{A_i(\zeta_2) - \mathrm{E}A_1(\zeta_2)\}\right|$$

$$\leq \sup_{|\zeta_1-\zeta_2|<\rho} \frac{1}{\sqrt{k}} \sum_{i=1}^{n} |A_i(\zeta_2) - A_i(\zeta_1)| + \sup_{|\zeta_1-\zeta_2|<\rho} \frac{n}{\sqrt{k}} \mathrm{E}|A_1(\zeta_2) - A_1(\zeta_1)|,$$

so that, due to Lemma 4,

$$\lim_{\rho \to 0} \limsup_{n} P\left(\sup_{|\zeta_1-\zeta_2|<\rho} \sup_{0 \le t \le 1} |M_n(t,\zeta_1) - M_n(t,\zeta_2)| > \varepsilon\right) = 0. \quad (4.14)$$

Then, since

$$P\left(\sup_{\zeta \in [-K,K]} \sup_{0 \le t \le 1} |M_n(t,\zeta) - M_n(t,0)| > \varepsilon\right)$$
$$\le P\left(\max_{j \in \mathbb{Z}} \sup_{0 \le t \le 1} |M_n(t,j\rho) - M_n(t,0)| > \frac{\varepsilon}{2}\right)$$
$$+ P\left(\sup_{|\zeta_1-\zeta_2|<\rho} \sup_{0 \le t \le 1} |M_n(t,\zeta_1) - M_n(t,\zeta_2)| > \frac{\varepsilon}{2}\right),$$

where the maximum is taken over $j$'s with $j\rho \in [-K,K]$, the lemma is asserted by (4.6) and (4.14). □

**Lemma 6.** *Suppose that conditions* **(A1)–(A3)** *hold. If* $\varphi(x) = \varphi_1(x) := I(x > 0)$, *then*

$$\frac{1}{\sqrt{1+\omega}} M_n(\cdot) \Rightarrow B(\cdot) \quad \text{in } D[0,1].$$

*If* $\varphi(x) = \varphi_2(x) := x_+$ *and (2.15) holds, then*

$$\frac{\alpha}{\sqrt{2+\chi}} M_n(\cdot) \Rightarrow B(\cdot) \quad \text{in } D[0,1].$$

**Proof.** We first prove the theorem for the case $\varphi(x) = I(x > 0)$. Let $m_n = [n/r_n]$. For $0 < \varepsilon < 1$, we set $I_j = I_{nj}(\varepsilon) = \{(j-1)r_n + 1, \ldots, (j-1)r_n + [(1-\varepsilon)r_n]\}$ and $J_j = J_{nj}(\varepsilon) = \{(j-1)r_n + 1, \ldots, jr_n\} \cap I_j^c$. Further, we set $m_n(t)$ as the integer part of $[nt]/r_n$. Then we can express

$$\frac{1}{\sqrt{1+\omega}} M_n(t,0) = \frac{1}{\sqrt{k(1+\omega)}} \left\{\sum_{j=1}^{m_n(t)} \sum_{i \in I_j} \{A_i(0) - \mathrm{E}A_i(0)\} + \sum_{j=1}^{m_n(t)} \sum_{i \in J_j} \{A_i(0) - \mathrm{E}A_i(0)\}\right\}$$
$$+ Z_n(t),$$

where

$$Z_n(t) = \frac{1}{\sqrt{k(1+\omega)}} \sum_{i=r_n m_n(t)+1}^{[nt]} \{A_i(0) - \mathrm{E}A_i(0)\}.$$

Owing to (2.8) and the fact that $|A_i(0)| \le 1$, we can see that

$$\|Z_n\| := \sup_{t \in [0,1]} |Z_n(t)| = o_P(1). \quad (4.15)$$

Set
$$L_{n,\varepsilon}(t) := \frac{1}{\sqrt{k(1+\omega)(1-\varepsilon)}} \sum_{j=1}^{m_n(t)} Y_j,$$

where $Y_j := \sum_{i \in I_j}\{A_i(0) - \mathrm{E}A_i(0)\}$ is $\{X_i : i \in I_j\}$-measurable. Let $D$ denote a space of cadlag functions on $[0,1]$ endowed with Skorokhod's metric (cf. Billingsley (1999)) and let $\mathcal{D}$ be its Borel $\sigma$-field. Define

$$\tilde{M}_n(t; x_1, x_2, \ldots, x_{m_n}) := \frac{1}{\sqrt{k(1+\omega)(1-\varepsilon)}} \sum_{i=1}^{m_n(t)} x_i.$$

Since $\tilde{M}_n : (\mathbb{R}^{m_n}, \mathcal{R}^{m_n}) \to (D, \mathcal{D})$ is measurable, owing to (2.7), we have that for $H \in \mathcal{D}$,

$$\begin{aligned} P(L_{n,\varepsilon} \in H) &= P((Y_1, \ldots, Y_{m_n}) \in \tilde{M}_n^{-1}(H)) \\ &= \tilde{P}((\tilde{Y}_1, \ldots, \tilde{Y}_{m_n}) \in \tilde{M}_n^{-1}(H)) + \mathrm{o}(1) \\ &= \tilde{P}(\tilde{L}_{n,\varepsilon} \in H) + \mathrm{o}(1), \end{aligned} \quad (4.16)$$

where
$$\tilde{L}_{n,\varepsilon}(t) := \frac{1}{\sqrt{k(1+\omega)(1-\varepsilon)}} \sum_{j=1}^{m_n(t)} \tilde{Y}_j,$$

and $\tilde{Y}_j, j = 1, \ldots, m_n$, are i.i.d. copies of $Y_1$ (cf. Eberlein (1984)). Since by (2.8) and (2.10),

$$\tilde{L}_{n,\varepsilon} \Rightarrow B \quad \text{in } D[0,1] \quad (4.17)$$

(cf. Theorem 18.2 of Billingsley (1999)), in view of (4.16), we have

$$L_{n,\varepsilon} \Rightarrow B \quad \text{in } D[0,1]. \quad (4.18)$$

Similarly, it can be verified that

$$N_{n,\varepsilon}(t) := \frac{1}{\sqrt{k(1+\omega)\varepsilon}} \sum_{j=1}^{m_n(t)} \sum_{i \in J_j} (A_i(0) - \mathrm{E}A_i(0)) \Rightarrow B \quad \text{in } D[0,1]. \quad (4.19)$$

Now, if we set $V_{n,\varepsilon} := (\sqrt{1-\varepsilon} - 1)L_{n,\varepsilon} + \sqrt{\varepsilon}N_{n,\varepsilon} + Z_n$, we can express

$$\frac{1}{\sqrt{1+\omega}} M_n = L_{n,\varepsilon} + V_{n,\varepsilon},$$

and thus, for any closed set $H \subset D$ and $\delta > 0$,

$$\limsup_n P\left(\frac{1}{\sqrt{1+\omega}} M_n \in H\right) \leq \limsup_n P(L_{n,\varepsilon} \in \bar{H}^\delta) + \limsup_n P(\|V_{n,\varepsilon}\| \geq \delta),$$

where $\|V_{n,\varepsilon}\| = \sup_{t\in[0,1]} |V_{n,\varepsilon}(t)|$ and $\bar{H}^\delta$ is the closure of $H^\delta$. From (4.18) and (4.19), we have

$$\lim_{\varepsilon \to 0} \limsup_n P((1 - \sqrt{1-\varepsilon})\|L_{n,\varepsilon}\| \geq \delta) = 0$$

and

$$\lim_{\varepsilon \to 0} \limsup_n P(\sqrt{\varepsilon}\|N_{n,\varepsilon}\| \geq \delta) = 0,$$

which, together with (4.15), yields that

$$\lim_{\varepsilon \to 0} \limsup_n P(\|V_{n,\varepsilon}\| \geq \delta) = 0.$$

Thus, by letting $\varepsilon \to 0$, we get

$$\limsup_n P\left(\frac{1}{\sqrt{1+\omega}} M_n \in H\right) \leq P(B \in \bar{H}^\delta).$$

Further, by letting $\delta \downarrow 0$, we have

$$\limsup_n P\left(\frac{1}{\sqrt{1+\omega}} M_n \in H\right) \leq P(B \in H),$$

which entails $(1+\omega)^{-1/2} M_n \Rightarrow B$ in $D[0,1]$ due to Portmanteau's theorem.

Next, we deal with the case that $\varphi(x) = x_+$. We first demonstrate that (4.15) still holds for this case. Note that

$$\left|\frac{1}{\sqrt{k}} \sum_{i=r_n m_n(t)+1}^{[nt]} \{A_i(0) - \mathrm{E}A_i(0)\}\right| \leq \frac{1}{\sqrt{k}} \sum_{i=r_n m_n(t)+1}^{[nt]} A_i(0) + \frac{1}{\sqrt{k}} \sum_{i=r_n m_n(t)+1}^{[nt]} \mathrm{E}A_i(0)$$

$$\leq \frac{1}{\sqrt{k}} \sum_{i=r_n m_n(t)+1}^{[nt]} A_i(0) + \frac{r_n}{\sqrt{k}}, \quad (4.20)$$

and

$$\sup_{t\in[0,1]} \frac{1}{\sqrt{k}} \sum_{i=r_n m_n(t)+1}^{[nt]} A_i(0) \leq \max_{1\leq j\leq m_n+1} \frac{1}{\sqrt{k}} \sum_{i=(j-1)r_n+1}^{jr_n} A_i(0). \quad (4.21)$$

We have that for every $\eta > 0$,

$$P\left(\max_{1\leq j\leq m_n+1} \frac{1}{\sqrt{k}} \sum_{i=(j-1)r_n+1}^{jr_n} A_i(0) > \eta\right) \leq (m_n+1) P\left(\frac{1}{\sqrt{k}} \sum_{i=1}^{r_n} A_i(0) > \eta\right),$$

and further, for sufficiently large $n$,

$$
\begin{aligned}
P\left(\frac{1}{\sqrt{k}} \sum_{i=1}^{r_n} (\log X_i - \log b(n/k))_+ > \eta\right) &\leq r_n P\left(\log X_1 - \log b(n/k) > \frac{\eta\sqrt{k}}{r_n}\right) \\
&= r_n \bar{F}(e^{\eta\sqrt{k}/r_n} b(n/k)) \\
&\sim \frac{r_n k}{n} e^{-\alpha\eta\sqrt{k}/r_n} \frac{l(e^{\eta\sqrt{k}/r_n} b(n/k))}{l(b(n/k))} \\
&\leq \frac{r_n}{n} \cdot Ck e^{-(\alpha-\delta)\eta\sqrt{k}/r_n}
\end{aligned}
\quad (4.22)
$$

for some $0 < \delta < \alpha$ and $C > 1$, where we have used Potter's theorem (cf. Bingham *et al.* (1987)). Thus, by (2.15), we get

$$
\max_{1 \leq j \leq m_n+1} \frac{1}{\sqrt{k}} \sum_{i=(j-1)r_n+1}^{jr_n} A_i(0) = o_P(1),
$$

which, together with (4.20) and (4.21), entails (4.15).

Next, we verify that

$$
\frac{\alpha}{\sqrt{k(2+\chi)(1-\varepsilon)}} \sum_{j=1}^{m_n(\cdot)} \sum_{i \in I_j} \{A_i(0) - \mathrm{E}A_i(0)\} \Rightarrow B(\cdot) \quad \text{in } D[0,1], \quad (4.23)
$$

which corresponds to (4.17) in the case of $\varphi(x) = I(x > 0)$. Let

$$
\xi_{nj} = \frac{\alpha}{\sqrt{k(2+\chi)(1-\varepsilon)}} \sum_{i \in I_j} \{A_i(0) - \mathrm{E}A_1(0)\} \quad \text{and} \quad \sigma_{nj}^2 = \mathrm{Var}(\xi_{nj}).
$$

Since $\frac{n}{k} \mathrm{E}(\log X_1 - \log b(n/k))_+^2 \sim 2\alpha^{-2}$, by (2.9), we have

$$
\sum_{j=1}^{[m_n t]} \sigma_{nj}^2 \longrightarrow t \quad \text{as } n \to \infty.
$$

Further, since $\frac{n}{k} \mathrm{E}(\log X_1 - \log b(n/k))_+^4 \sim 4!\alpha^{-4}$, we have that for some $K_1 > 0$,

$$
m_n \mathrm{E}(\xi_{n1}^4) \leq \frac{n}{r_n} K_1 \frac{r_n^4}{k^2} \frac{k}{n} \leq K_1 \frac{r_n^3}{k}.
$$

Henceforth, by using (4.22) and the Schwarz inequality, for each $\eta > 0$, there exist some $c, K > 0$, such that

$$
\sum_{j=1}^{m_n} \mathrm{E}[\xi_{nj}^2 I(|\xi_{nj}| \geq \eta)] = m_n \mathrm{E}[\xi_{n1}^2 I(|\xi_{n1}| \geq \eta)] \leq m_n \sqrt{\mathrm{E}(\xi_{n1}^4) P(|\xi_{n1}| \geq \eta)}
$$

$$= \sqrt{m_n \mathrm{E}(\xi_{n1}^4) m_n P(|\xi_{n1}| \geq \eta)} \leq \sqrt{K r_n^3 \mathrm{e}^{-c\sqrt{k}/r_n}},$$

where the last term is o(1) by (2.15). By the same reasoning as in the derivation of (4.16), we can view $\{\xi_{nj} : j = 1, \ldots, m_n, n = 1, 2, \ldots\}$ as a row-wise independent double array of zero-mean r.v.'s. Therefore, by Theorem 18.2 of Billingsley (1999), we can obtain (4.23). Since the rest of the proof essentially follows the same lines below (4.17), we omit it for brevity. □

**Proof of Theorem 1.** According to Lemma 6, (2.14) and (2.16) are asserted if we verify that

$$\frac{1}{\sqrt{1+\omega}} T_n(\varphi_1) = \frac{1}{\sqrt{k(1+\omega)}} \max_{1 \leq l \leq n} \left| \sum_{i=1}^{l} I(X_i > b(n/k)) - \frac{l}{n} \sum_{i=1}^{n} I(X_i > b(n/k)) \right| \\ + \mathrm{o}_P(1), \quad (4.24)$$

and

$$\frac{\alpha}{\sqrt{2+\chi}} T_n(\varphi_2) = \frac{\alpha}{\sqrt{k(2+\chi)}} \max_{1 \leq l \leq n} \left| \sum_{i=1}^{l} (\log X_i - \log b(n/k))_+ \right. \\ \left. - \frac{l}{n} \sum_{i=1}^{n} (\log X_i - \log b(n/k))_+ \right| + \mathrm{o}_P(1). \quad (4.25)$$

Since the proof of (4.25) is similar to that of (4.24), we only provide for the latter. Let $K$ be any positive real number. By setting $\zeta_n = -\sqrt{k}(\log X_{(k)} - \log b(n/k))$, we can express

$$\frac{1}{\sqrt{k}} \max_{1 \leq l \leq n} \left| \sum_{i=1}^{l} \{I(X_i > X_{(k)}) - I(X_i > b(n/k))\} - \frac{l}{n} \sum_{i=1}^{n} \{I(X_i > X_{(k)}) - I(X_i > b(n/k))\} \right|$$

$$= \frac{1}{\sqrt{k}} \max_{1 \leq l \leq n} \left| \sum_{i=1}^{l} \{A_i(\zeta_n) - A_i(0)\} - \frac{l}{n} \sum_{i=1}^{n} \{A_i(\zeta_n) - A_i(0)\} \right|$$

$$:= I_n + II_n,$$

where

$$I_n = \frac{1}{\sqrt{k}} \max_{1 \leq l \leq n} \left| \sum_{i=1}^{l} \{A_i(\zeta_n) - A_i(0)\} - \frac{l}{n} \sum_{i=1}^{n} \{A_i(\zeta_n) - A_i(0)\} \right| I(|\zeta_n| < K)$$

and

$$II_n = \frac{1}{\sqrt{k}} \max_{1 \leq l \leq n} \left| \sum_{i=1}^{l} \{A_i(\zeta_n) - A_i(0)\} - \frac{l}{n} \sum_{i=1}^{n} \{A_i(\zeta_n) - A_i(0)\} \right| I(|\zeta_n| \geq K).$$

Due to Lemma 5, we have

$$I_n \leq \frac{1}{\sqrt{k}} \max_{1 \leq l \leq n} \sup_{\zeta \in [-K,K]} \left| \sum_{i=1}^{l} \{A_i(\zeta) - A_i(0)\} - \frac{l}{n} \sum_{i=1}^{n} \{A_i(\zeta) - A_i(0)\} \right|$$

$$\leq 2 \sup_{\zeta \in [-K,K]} \sup_{0 \leq t \leq 1} |M_n(t,\zeta) - M_n(t,0)| = o_P(1).$$

On the other hand, since $\zeta_n \Rightarrow N(0, \alpha^{-2}(1+\omega))$ in view of Lemmas 3 and 6 (cf. Theorem 2.4 of Hsing (1991)), we have

$$\limsup_n P(II_n > 0) \leq \limsup_n P(|\zeta_n| \geq K) \to 0 \qquad \text{as } K \to \infty.$$

Therefore, (4.24) is asserted. This completes the proof. □

Below, we prove Theorem 2. It is well known that under conditions **(B1)**–**(B3)**, the sequence of stochastic processes $\mathcal{E}_n$ defined by

$$\mathcal{E}_n(x) = \sqrt{k} \left( \frac{1}{k} \sum_{i=1}^{n} I(Z_i > x^{-1/\alpha} b^*(n/k)) - x \right), \qquad x \in [0, \infty), \tag{4.26}$$

converges weakly to a standard Brownian motion $B$ in $D[0, \infty)$ (cf. Proposition 2.1 of Resnick *et al.* (1997b)). The following result is due to Proposition 3.2 of Resnick *et al.* (1997a), which plays an important role in verifying Theorem 2.

**Lemma 7.** *Under conditions* **(B1)**–**(B4)**,

$$\frac{1}{\sqrt{k}} \sup_{x \in [c,d]} \sum_{i=1}^{n} |I(Z_i > xb^*(n/k)) - I(\hat{Z}_i > xb^*(n/k))| = o_P(1)$$

*for every* $0 < c < d < \infty$.

**Lemma 8.** *Under conditions* **(B1)**–**(B4)**,

$$\sqrt{k} \{\log \hat{Z}_{(k)} - \log b^*(n/k)\} = O_P(1).$$

**Proof.** By Lemmas 3, 6 and 7, for every $\zeta \in \mathbb{R}$,

$$\frac{1}{\sqrt{k}} \sum_{i=1}^{n} \{I(\hat{Z}_i > e^{\zeta/\sqrt{k}} b^*(n/k)) - P(Z_i > e^{\zeta/\sqrt{k}} b^*(n/k))\}$$

$$= \frac{1}{\sqrt{k}} \sum_{i=1}^{n} \{I(Z_i > e^{\zeta/\sqrt{k}} b^*(n/k)) - P(Z_i > e^{\zeta/\sqrt{k}} b^*(n/k))\} + o_P(1) \Rightarrow N(0,1).$$

Hence,
$$\sqrt{k}\{\log \hat{Z}_{(k)} - \log b^*(n/k)\} \Rightarrow N(0, \alpha^{-2})$$
(cf. Theorem 2.4 of Hsing (1991)). This completes the proof. □

**Lemma 9.** *Under conditions* **(B1)**–**(B4)** *and (2.23), for $c < d \in \mathbb{R}$,*

$$\frac{1}{\sqrt{k}} \sup_{x \in [c,d]} \sum_{i=1}^{n} |(\log Z_i - \log b^*(n/k) + x)_+ - (\log \hat{Z}_i - \log b^*(n/k) + x)_+| = o_P(1). \quad (4.27)$$

**Proof.** By setting $Y_i = \log Z_i - \log b^*(n/k)$ and $\hat{Y}_i = \log \hat{Z}_i - \log b^*(n/k)$, we can express the argument in (4.27) as

$$\begin{aligned}
&\frac{1}{\sqrt{k}} \sup_{x \in [c,d]} \sum_{i=1}^{n} |(Y_i + x)_+ - (\hat{Y}_i + x)_+| \\
&\leq \frac{1}{\sqrt{k}} \sup_{x \in [c,d]} \sum_{i=1}^{n} |(Y_i + x)_+ - (\hat{Y}_i + x)_+||I(\hat{Y}_i + x > 0) - I(Y_i + x > 0)| \\
&\quad + \frac{1}{\sqrt{k}} \sup_{x \in [c,d]} \sum_{i=1}^{n} |(Y_i + x)_+ - (\hat{Y}_i + x)_+|I(\hat{Y}_i + x > 0, Y_i + x > 0) \\
&:= \Lambda_{n1} + \Lambda_{n2}.
\end{aligned} \quad (4.28)$$

Note that $\Lambda_{n1}$ is no more than $I_n + II_n$, where

$$I_n = \frac{1}{\sqrt{k}} \sup_{x \in [c,d]} \sum_{i=1}^{n} U_i(x) I((Y_i + x)_+ \vee (\hat{Y}_i + x)_+ \leq 1),$$

$$II_n = \frac{1}{\sqrt{k}} \sup_{x \in [c,d]} \sum_{i=1}^{n} U_i(x) \{I((Y_i + x)_+ > 1, (\hat{Y}_i + x)_+ \leq 0)$$
$$+ I((Y_i + x)_+ \leq 0, (\hat{Y}_i + x)_+ > 1)\}$$

and

$$U_i(x) = |(Y_i + x)_+ - (\hat{Y}_i + x)_+||I(\hat{Y}_i + x > 0) - I(Y_i + x > 0)|.$$

First, note that by Lemma 7,

$$I_n \leq \frac{1}{\sqrt{k}} \sup_{x \in [c,d]} \sum_{i=1}^{n} |I(\hat{Y}_i + x > 0) - I(Y_i + x > 0)| = o_P(1).$$

Second, if we set $\mathbf{X}_{i-1} = (X_{i-1}, \ldots, X_{i-p})^T$ and $|\mathbf{x}| = \sqrt{\mathbf{x}^T\mathbf{x}}$ for $\mathbf{x} \in \mathbb{R}^p$, we have that for $0 < \delta < \nu$,

$$\frac{1}{n^{1/\alpha+\delta}} \max_{1 \leq i \leq n} |\mathbf{X}_{i-1}| = o_P(1),$$

and therefore, by using (2.23) and the fact that

$$b^*(\cdot) \in RV_{1/\alpha} \tag{4.29}$$

(cf. Theorems 1.5.12 and 1.5.4 of Bingham *et al.* (1987)), we get

$$\sup_{x \in [c,d]} \max_{1 \leq i \leq n} \{I((Y_i+x)_+ > 1, (\hat{Y}_i+x)_+ \leq 0) + I((Y_i+x)_+ \leq 0, (\hat{Y}_i+x)_+ > 1)\}$$

$$\leq I\left(\max_{1 \leq i \leq n} |Z_i - \hat{Z}_i| \geq e^{-d}b^*(n/k)(e-1)\right)$$

$$\leq I\left(\frac{n^{1/\alpha+\delta}}{d(n)b^*(n/k)} \cdot d(n)|\hat{\phi} - \phi| \frac{1}{n^{1/\alpha+\delta}} \max_{1 \leq i \leq n} |\mathbf{X}_{i-1}| > e^{-d}(e-1)\right)$$

$$= o_P(1),$$

which asserts $II_n = o_P(1)$. Hence $\Lambda_{n1} = o_P(1)$.

Third, by using (2.23), (4.29), Theorem 1.5.4 of Bingham *et al.* (1987) and the fact that

$$|(\log z - \log b^*(n/k) + x)_+ - (\log \hat{z} - \log b^*(n/k) + x)_+| \leq \frac{|z - \hat{z}|}{\min\{z, \hat{z}\}},$$

we have

$$\frac{1}{\sqrt{k}} \sup_{x \in [c,d]} \sum_{i=1}^{n} |(Y_i+x)_+ - (\hat{Y}_i+x)_+| I(\hat{Y}_i + x > 0, Y_i + x > 0)$$

$$\leq \frac{\max_{1 \leq i \leq n} |Z_i - \hat{Z}_i|}{e^{-d}b^*(n/k)} \frac{1}{\sqrt{k}} \sum_{i=1}^{n} I(Z_i > e^{-d}b^*(n/k))$$

$$\leq e^d \frac{\sqrt{k} n^{1/\alpha+\delta}}{d(n)b^*(n/k)} d(n)|\hat{\phi} - \phi| \frac{1}{n^{1/\alpha+\delta}} \max_{1 \leq i \leq n} |\mathbf{X}_{i-1}| \frac{1}{k} \sum_{i=1}^{n} I(Z_i > e^{-d}b^*(n/k)) \tag{4.30}$$

$$= o_P(1).$$

This implies $\Lambda_{n2} = o_P(1)$. Hence the lemma is established by (4.28). □

**Proof of Theorem 2.** Note that

$$\frac{1}{\sqrt{k}} \max_{1 \leq l \leq n} \left| \sum_{i=1}^{l} \{I(\hat{Z}_i > \hat{Z}_{(k)}) - I(Z_i > b^*(n/k))\} \right|$$

$$-\frac{l}{n}\sum_{i=1}^{n}\{I(\hat{Z}_i > \hat{Z}_{(k)}) - I(Z_i > b^*(n/k))\}\bigg|$$

is bounded by $I_n + II_n$, where

$$I_n = \frac{1}{\sqrt{k}} \max_{1 \le l \le n} \bigg|\sum_{i=1}^{l}\{I(\hat{Z}_i > \hat{Z}_{(k)}) - I(Z_i > \hat{Z}_{(k)})\} - \frac{l}{n}\sum_{i=1}^{n}\{I(\hat{Z}_i > \hat{Z}_{(k)}) - I(Z_i > \hat{Z}_{(k)})\}\bigg|,$$

and

$$II_n = \frac{1}{\sqrt{k}} \max_{1 \le l \le n} \bigg|\sum_{i=1}^{l}\{I(Z_i > \hat{Z}_{(k)}) - I(Z_i > b^*(n/k))\}$$
$$- \frac{l}{n}\sum_{i=1}^{n}\{I(Z_i > \hat{Z}_{(k)}) - I(Z_i > b^*(n/k))\}\bigg|.$$

By Lemmas 7 and 8, we can have

$$I_n \le \frac{2}{\sqrt{k}}\sum_{i=1}^{n}\bigg|I\bigg(\hat{Z}_i > \frac{\hat{Z}_{(k)}}{b^*(n/k)}b^*(n/k)\bigg) - I\bigg(Z_i > \frac{\hat{Z}_{(k)}}{b^*(n/k)}b^*(n/k)\bigg)\bigg| = o_P(1).$$

Further, by using Lemma 8, we can prove $II_n = o_P(1)$ in a fashion similar to that used to prove Theorem 1. Therefore, in view of Corollary 1, we have

$$T_n^*(\varphi_1) = \frac{1}{\sqrt{k}} \max_{1 \le l \le n} \bigg|\sum_{i=1}^{l}I(Z_i > b^*(n/k)) - \frac{l}{n}\sum_{i=1}^{n}I(Z_i > b^*(n/k))\bigg| + o_P(1)$$
$$\Rightarrow \sup_{0 \le t \le 1} |B^\circ(t)|,$$

which establishes (2.22). Since (2.24) can be proven similarly – in this case Lemma 9 is used instead of Lemma 7 – we complete the proof without detailing algebras. □

Now we prove Theorem 3. It is not difficult to verify the following lemma (cf. Theorem 3.1 of Hsing (1991)).

**Lemma 10.** *Under the conditions in Theorem 3, we have that for each $t \in (0, \tau]$ and $\delta > 0$,*

$$\sup_{x > \delta}\bigg|\frac{1}{k}\sum_{i=1}^{[nt]}I(X_i > xb_{\text{pre}}(n/k)) - tx^{-\alpha_{\text{pre}}}\bigg| \xrightarrow{P} 0$$

*and*

$$\sup_{x > \delta}\bigg|\frac{1}{k}\sum_{i=1}^{[nt]}(\log X_i - \log b_{\text{pre}}(n/k) - \log x)_+ - tx^{-\alpha_{\text{pre}}}\alpha_{\text{pre}}^{-1}\bigg| \xrightarrow{P} 0.$$

*Further, for each $t \in (\tau, 1]$,*

$$\sup_{x>\delta}\left|\frac{1}{k}\sum_{i=[n\tau]+1}^{[nt]} I(X_i > xb_{\text{post}}(n/k)) - (t-\tau)x^{-\alpha_{\text{post}}}\right| \xrightarrow{P} 0$$

*and*

$$\sup_{x>\delta}\left|\frac{1}{k}\sum_{i=[n\tau]+1}^{[nt]} (\log X_i - \log b_{\text{post}}(n/k) - \log x)_+ - (t-\tau)x^{-\alpha_{\text{post}}}\alpha_{\text{post}}^{-1}\right| \xrightarrow{P} 0.$$

**Proof of Theorem 3.** We verify the first part of the theorem. Suppose that (2.26) holds. Let $y > 1$ be a real number such that

$$\tau y^{-\alpha_{\text{pre}}} + \left(\frac{y}{c}\right)^{-\alpha_{\text{post}}}(1-\tau) > 1.$$

For sufficiently large $n$, we have

$$\frac{1}{k}\sum_{i=1}^{n} I(X_i > yb_{\text{pre}}(n/k)) = \frac{1}{k}\sum_{i=1}^{[n\tau]} I(X_i > yb_{\text{pre}}(n/k)) + \frac{1}{k}\sum_{i=[n\tau]+1}^{n} I(X_i > yb_{\text{pre}}(n/k))$$

$$\geq \frac{1}{k}\sum_{i=1}^{[n\tau]} I(X_i > yb_{\text{pre}}(n/k)) + \frac{1}{k}\sum_{i=[n\tau]+1}^{n} I\left(X_i > \frac{y}{c}b_{\text{post}}(n/k)\right)$$

$$\xrightarrow{P} \tau y^{-\alpha_{\text{pre}}} + \left(\frac{y}{c}\right)^{-\alpha_{\text{post}}}(1-\tau) > 1,$$

which implies $\lim_{n\to\infty} P(X_{(k)} > yb_{\text{pre}}(n/k)) = 1$. Hence, if we set

$$U_{n,1}(t) = \frac{1}{k}\sum_{i=1}^{[nt]} I(X_i > X_{(k)}) - \frac{[nt]}{n}\left(1 - \frac{1}{k}\right),$$

we have

$$U_{n,1}(\tau) \leq \frac{1}{k}\sum_{i=1}^{[n\tau]} I(X_i > yb_{\text{pre}}(n/k)) - \tau + o_P(1) = (y^{-\alpha_{\text{pre}}} - 1)\tau + o_P(1).$$

Since $T_n(\varphi_1) = \sup_{0 \leq t \leq 1} \sqrt{k}|U_{n,1}(t)|$ and $(y^{-1/\alpha_{\text{pre}}} - 1) < 0$, the above asserts that $T_n(\varphi_1) \xrightarrow{P} \infty$.

Next, we verify the second part of the theorem. We first handle the case that $d$ in (2.27) is finite. By Lemma 10, we have that for $0 < y < \infty$,

$$\frac{1}{k}\sum_{i=1}^{n} I(X_i > y b_{\text{pre}}(n/k)) = \frac{1}{k}\sum_{i=1}^{[n\tau]} I(X_i > y b_{\text{pre}}(n/k))$$

$$+ \frac{1}{k}\sum_{i=[n\tau]+1}^{n} I\left(X_i > y \frac{b_{\text{pre}}(n/k)}{b_{\text{post}}(n/k)} b_{\text{post}}(n/k)\right)$$

$$\xrightarrow{P} \tau y^{-\alpha_{\text{pre}}} + \left(\frac{y}{d}\right)^{-\alpha_{\text{post}}}(1-\tau),$$

and thus there exists $y_0 > 1$ such that $\tau y_0^{-\alpha_{\text{pre}}} + (\frac{y_0}{d})^{-\alpha_{\text{post}}}(1-\tau) = 1$,

$$\frac{X_{(k)}}{b_{\text{pre}}(n/k)} \xrightarrow{P} y_0 \quad \text{and} \quad \frac{X_{(k)}}{b_{\text{post}}(n/k)} \xrightarrow{P} \frac{y_0}{d}.$$

Hence, we have that for each $t \in (0, \tau]$,

$$\frac{1}{k}\sum_{i=1}^{[nt]} I(X_i > X_{(k)}) = \frac{1}{k}\sum_{i=1}^{[nt]} I\left(X_i > \frac{X_{(k)}}{b_{\text{pre}}(n/k)} b_{\text{pre}}(n/k)\right) \xrightarrow{P} t y_0^{-\alpha_{\text{pre}}}$$

and

$$\frac{1}{k}\sum_{i=1}^{[nt]} (\log X_i - \log X_{(k)})_+ = \frac{1}{k}\sum_{i=1}^{[nt]} \left(\log X_i - \log b_{\text{pre}}(n/k) - \log \frac{X_{(k)}}{b_{\text{pre}}(n/k)}\right)_+ \xrightarrow{P} \frac{t y_0^{-\alpha_{\text{pre}}}}{\alpha_{\text{pre}}}.$$

Further, for each $t \in (\tau, 1]$,

$$\frac{1}{k}\sum_{i=[n\tau]+1}^{[nt]} I(X_i > X_{(k)}) = \frac{1}{k}\sum_{i=[n\tau]+1}^{[nt]} I\left(X_i > \frac{X_{(k)}}{b_{\text{post}}(n/k)} b_{\text{post}}(n/k)\right) \xrightarrow{P} (t-\tau)\left(\frac{y_0}{d}\right)^{-\alpha_{\text{post}}}$$

and

$$\frac{1}{k}\sum_{i=[n\tau]+1}^{[nt]} (\log X_i - \log X_{(k)})_+ = \frac{1}{k}\sum_{i=[n\tau]+1}^{[nt]} \left(\log X_i - \log b_{\text{post}}(n/k) - \log \frac{X_{(k)}}{b_{\text{post}}(n/k)}\right)_+$$

$$\xrightarrow{P} (t-\tau)\left(\frac{y_0}{d}\right)^{-\alpha_{\text{post}}} \frac{1}{\alpha_{\text{post}}}.$$

due to Lemma 10. Since both the stochastic process $\frac{1}{k}\sum_{i=1}^{[nt]} I(X_i > X_{(k)})$ and its limiting process are non-decreasing in $t$ and the limiting process is continuous in $t$, we have

$$\sup_{t\in[0,\tau]}\left|\frac{1}{k}\sum_{i=1}^{[nt]} I(X_i > X_{(k)}) - ty_0^{-\alpha_{\text{pre}}}\right| \xrightarrow{P} 0.$$

Similarly, we get

$$\sup_{t\in[0,\tau]}\left|\frac{1}{k}\sum_{i=1}^{[nt]} (\log X_i - \log X_{(k)})_+ - \frac{ty_0^{-\alpha_{\text{pre}}}}{\alpha_{\text{pre}}}\right| \xrightarrow{P} 0,$$

$$\sup_{t\in(\tau,1]}\left|\frac{1}{k}\sum_{i=[n\tau]+1}^{[nt]} I(X_i > X_{(k)}) - (t-\tau)\left(\frac{y_0}{d}\right)^{-\alpha_{\text{post}}}\right| \xrightarrow{P} 0$$

and

$$\sup_{t\in(\tau,1]}\left|\frac{1}{k}\sum_{i=[n\tau]+1}^{[nt]} (\log X_i - \log X_{(k)})_+ - (t-\tau)\left(\frac{y_0}{d}\right)^{-\alpha_{\text{post}}}\frac{1}{\alpha_{\text{post}}}\right| \xrightarrow{P} 0.$$

By using these facts, we can have that uniformly in $t \in [0,1]$,

$$U_{n,1}(t) \xrightarrow{P} t(y_0^{-\alpha_{\text{pre}}} - 1) \vee \left[\tau\left\{y_0^{-\alpha_{\text{pre}}} - \left(\frac{y_0}{d}\right)^{-\alpha_{\text{post}}}\right\} + t\left\{\left(\frac{y_0}{d}\right)^{-\alpha_{\text{post}}} - 1\right\}\right],$$

and further,

$$U_{n,2}(t) := \frac{1}{k}\sum_{i=1}^{[nt]}(\log X_i - \log X_{(k)})_+ - \frac{[nt]}{nk}\sum_{i=1}^{n}(\log X_i - \log X_{(k)})_+$$

$$\xrightarrow{P} (t(1-\tau) \wedge (1-t)\tau)\left\{\frac{y_0^{-\alpha_{\text{pre}}}}{\alpha_{\text{pre}}} - \left(\frac{y_0}{d}\right)^{-\alpha_{\text{post}}}\frac{1}{\alpha_{\text{post}}}\right\},$$

which has a minimum at $t = \tau$. Hence, $T_n(\varphi_2) \xrightarrow{P} \infty$ and (2.28) is established.

Now, we deal with the case that $d = \infty$. Due to Lemma 10, we have that for $y > 0$,

$$\frac{1}{k}\sum_{i=1}^{n} I(X_i > yb_{\text{post}}(n/k)) \xrightarrow{P} (1-\tau)y^{-\alpha_{\text{post}}},$$

and subsequently,

$$\frac{X_{(k)}}{b_{\text{post}}(n/k)} \xrightarrow{P} (1-\tau)^{1/\alpha_{\text{post}}} \quad \text{and} \quad \frac{X_{(k)}}{b_{\text{pre}}(n/k)} \xrightarrow{P} \infty.$$

By using this and Lemma 10, it can be shown that

$$\sup_{t\in(0,\tau]} \frac{1}{k} \sum_{i=1}^{[nt]} I(X_i > X_{(k)}) \xrightarrow{P} 0,$$

$$\sup_{t\in(0,\tau]} \frac{1}{k} \sum_{i=1}^{[nt]} (\log X_i - \log X_{(k)})_+ \xrightarrow{P} 0,$$

$$\sup_{t\in(\tau,1]} \left| \frac{1}{k} \sum_{i=[n\tau]+1}^{[nt]} I(X_i > X_{(k)}) - \frac{t-\tau}{1-\tau} \right| \xrightarrow{P} 0$$

and

$$\sup_{t\in(\tau,1]} \left| \frac{1}{k} \sum_{i=[n\tau]+1}^{[nt]} (\log X_i - \log X_{(k)})_+ - \frac{t-\tau}{1-\tau} \frac{1}{\alpha_{\text{post}}} \right| \xrightarrow{P} 0.$$

Further, by using the above arguments, it can be obtained that

$$\sup_{t\in[0,1]} \left| U_{n,1}(t) - (-t) \vee \left( \frac{t-\tau}{1-\tau} - t \right) \right| \xrightarrow{P} 0$$

and

$$\sup_{t\in[0,1]} \left| U_{n,2}(t) - (-t) \vee \left( \frac{t-\tau}{1-\tau} - t \right) \frac{1}{\alpha_{\text{post}}} \right| \xrightarrow{P} 0.$$

Since $t \mapsto (-t) \vee (\frac{t-\tau}{1-\tau} - t)$ has a minimum at $t = \tau$, the theorem is established. $\square$

## Acknowledgements

We thank the Editor-in-Chief and the two anonymous referees for their valuable comments to improve the quality of the paper. This work was supported by grant No. R01-2006-000-10545-0 from the Basic Research Program of the Korea Science & Engineering Foundation.